\documentclass[11pt,oneside]{article}
\usepackage{amsfonts,amssymb,amsmath,verbatim,authblk,hyperref,graphicx,amsbsy,mathrsfs,caption,subcaption}

\newcommand{\beq}{\begin{equation}\ }
\newcommand{\eeq}{\end{equation}\ }

\newtheorem{define}{Definition}[section]

\newtheorem{example}{Example}[section]

\begin{document}

\author[1,2]{Demetris T. Christopoulos}
\affil[1]{\small{National and Kapodistrian University of Athens, Department of Economics}}
\affil[2]{\tt{dchristop@econ.uoa.gr}, \tt{dem.christop@gmail.com} }

\title{{\itshape On the effective radius of convergence for a given truncated power series expansion}}

\maketitle

\begin{abstract}
An effective radius of convergence is defined and computed for any truncated Taylor series. Applications to well known series are performed and is shown that a range of good coincidence for actual and approximative plot can always be found. For sufficient large degree of approximation the effective radius is also an estimation of the true non-infinite radius of convergence.\medskip
\end{abstract}


\section{Statement of the problem}
\begin{figure}
\begin{center}
\caption{$sin(x)$ and its $11^{th}$ degree Taylor polynomial}\label{fig:00}
\vspace{0.5in}  \includegraphics[width=6cm,height=6cm]{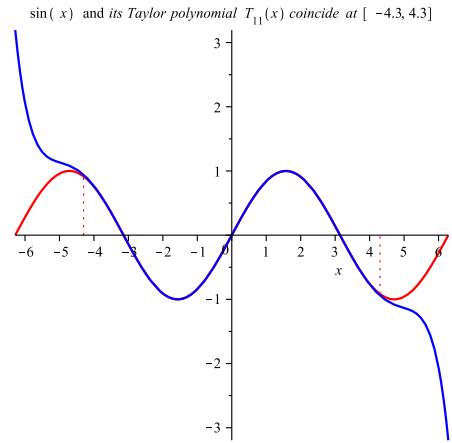}
\vspace{0.5in}
\end{center}
\end{figure}

 Let' s take the common $sin(x)$ function, its $11^{th}$ degree Taylor polynomial and plot them together, see Figure \ref{fig:00}. We directly observe that the two graphs coincide only for the relatively small interval $[-4.3,4.3]$, if we compare it with the $sinus$ infinite radius of convergence. 
So, when we have a series approximation of a known function around $x_0$, it is useful to be able to compute an interval with centre $x_0$ such that the function and its approximation practically coincide each other.  

\section{The truncated power series}
Many times we have only a truncated part of a general power series, either because the underlying function is unknown or because we simply do not want to take many terms of the known expansion. The question that arises naturally is `how far away' from the central point $x_{0}$ can we evaluate the truncated series with acceptable accuracy? Or in other words, given a level of accuracy, which is the radius with centre $x_0$ that will make the graphs of the real and the approximate function to `coincide'?. 
In order to be more rigorous we have to define the graph coincidence by use of a distance based on a proper norm.
\begin{define}\label{def:ecoinc}
Let two real functions $f,g$ defined in $[a,b]$ and let their graphs $G(f),G(g)$ for the same interval by using a proper $N$-points partition of it. Let $V_{f},\,V_{g}$ be the vectors for the values of $f,g$ at the partition and let $l$ a vector norm. We say that the graphs $(\epsilon,l)$-coincide at $[a,b]$ if and only if next inequality holds
\begin{equation}
\lVert V_{f}-V_{g} \rVert_{l}<\epsilon
\end{equation}
\end{define}
Since we know that a truncated series approximation has its biggest error at the highest distances from centre $x_0$, i.e. at $a$ or $b$ it is reasonable to choose the infinity vector norm for computing the distance between the two graphs.
Now we can give a formal definition of the $(\epsilon,l)$-effective radius of convergence.
\begin{define}
Let a truncated power series around $x_{0}$ be given as the sum
\begin{equation}
\label{eq:pn}
P_{m}(x)=\sum_{n=0}^{m}{a_{n}\,\left(x-x_{0}\right)^{n}}
\end{equation}
which approximates a function $f\in\,C^{\infty}\left[a,b\right]$. Then the $(\epsilon,l)$-effective radius of convergence $R_{ef}\left(\epsilon,l\right)$ is such that $\forall\,x\in\left(x_{0}-R_{ef},x_{0}+R_{ef} \right)$ the graphs $G(f)$, $G(P)$ $(\epsilon,l)$-coincide.
\end{define}

\begin{figure}
\caption{$P_{11}(x)$ of $sin(x)$ approximations}\label{fig:01} 
        \centering
        \begin{subfigure}[b]{0.4\textwidth}
                \includegraphics[width=\textwidth]{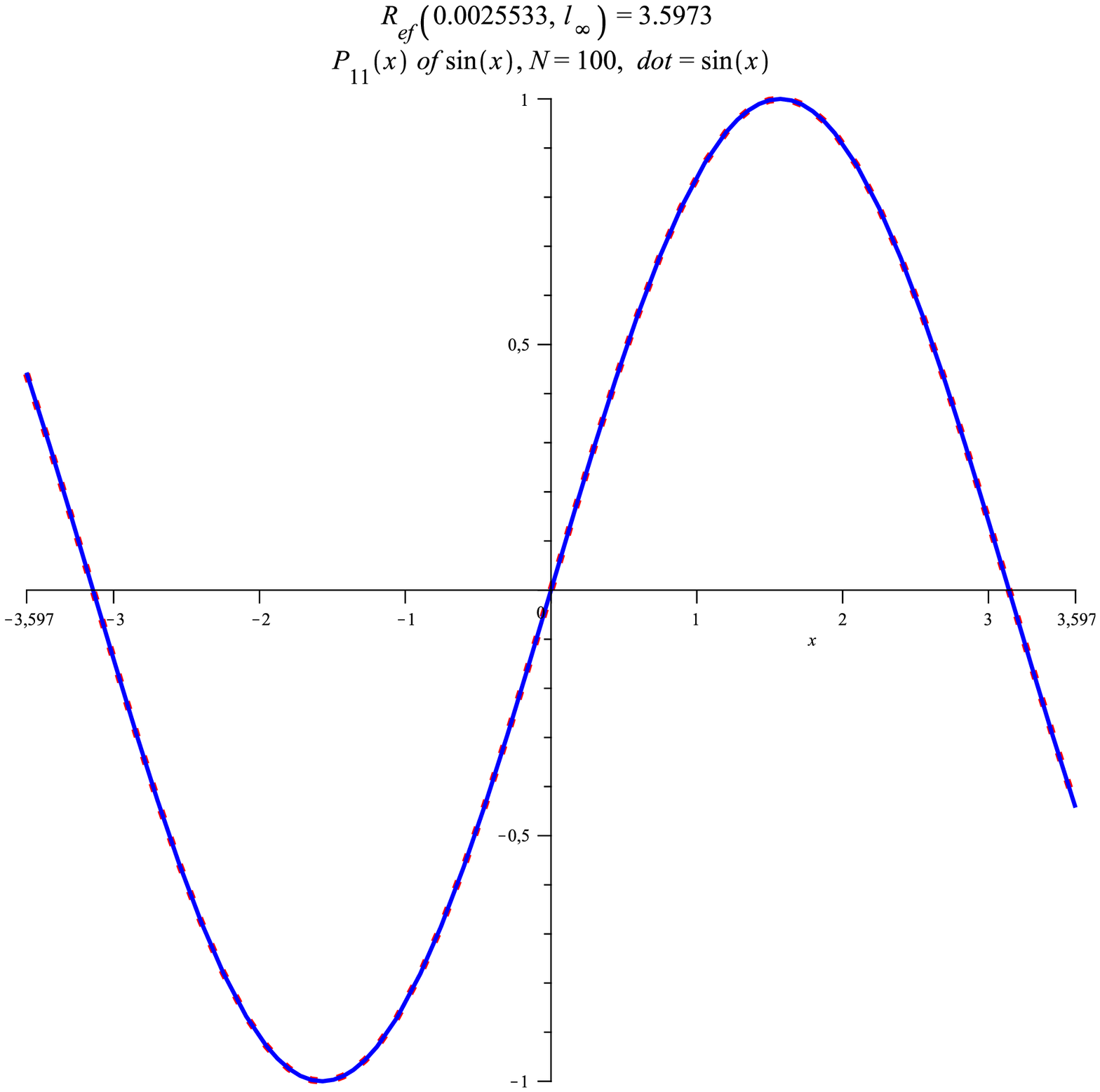}
                \caption{$R_{ef}(2.55\times10^{-3},l_{\infty})=3.5973$}
                \label{fig:01a}
        \end{subfigure}%
        \qquad
        \begin{subfigure}[b]{0.4\textwidth}
                \includegraphics[width=\textwidth]{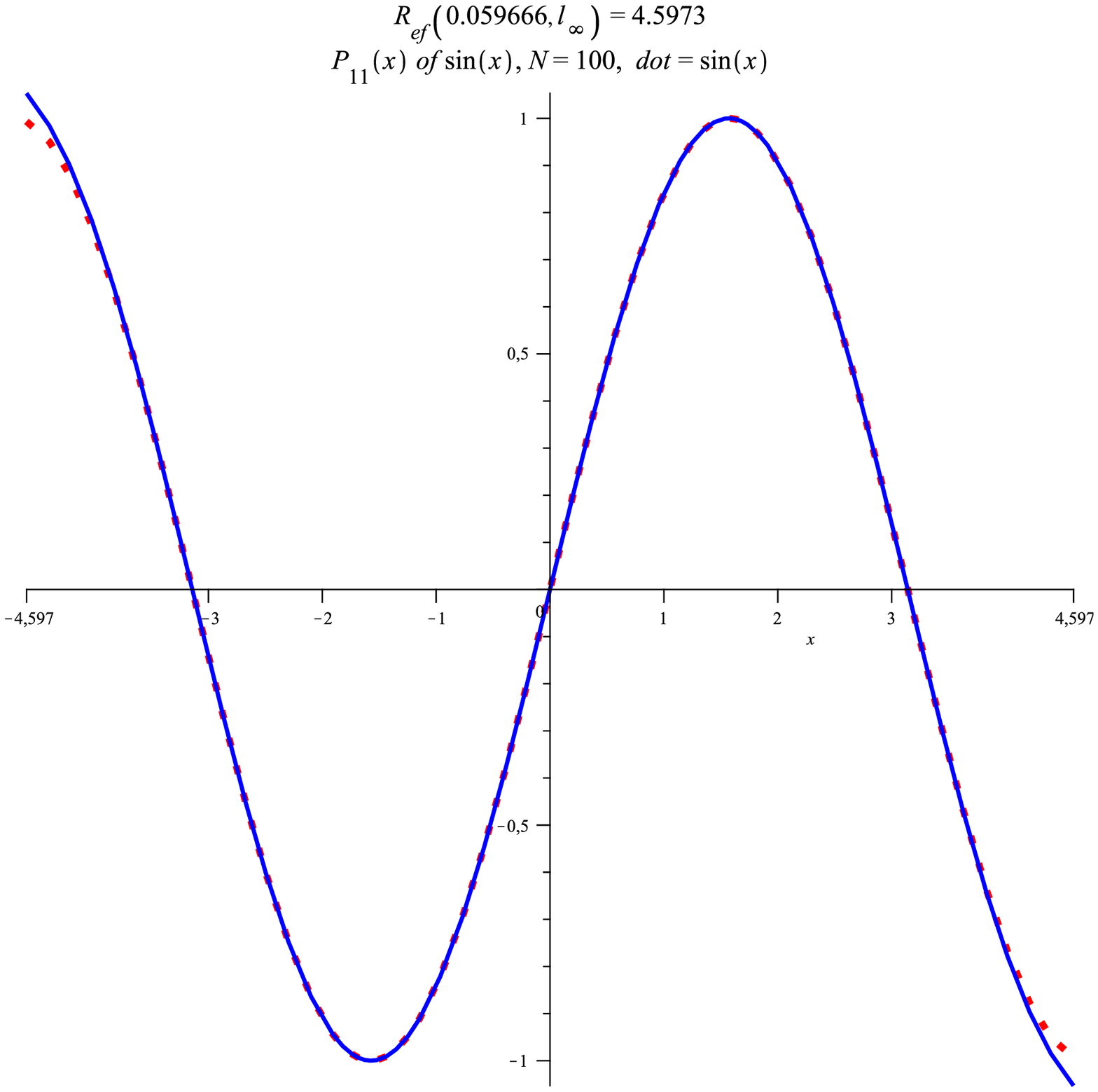}
                \caption{$R_{ef}(5.97\times10^{-2},l_{\infty})=4.5973$}
                \label{fig:01b}
        \end{subfigure}
      
\end{figure}

\begin{example}
For the MacLaurin power series of sinus we have that
\begin{equation}
\label{eq:sin}
sin(x)\approx \sum_{n=1}^{m}{(-1)^{n-1}\,\frac{x^{2n-1}}{\left(2n-1\right)!}}=
x-\frac{x^3}{3!}+\frac{x^5}{5!}-\frac{x^7}{7!}+\ldots+(-1)^{m}\frac{x^{2m-1}}{(2m-1)!}
\end{equation}
if we take $m=6$ we have
\begin{equation}
\label{eq:sin11}
P_{11}(x)= x-\frac{1}{6}\,{x}^{3}+{\frac {1}{120}}\,{x}^{5}-{
\frac {1}{5040}}\,{x}^{7}+{\frac {1}{362880}}\,{x}^{9}-{\frac {1}{39916800}}\,{x}^{11}
\end{equation}
The graphs of $sin(x),P_{11}(x)$, if we take $N=100$ points at $[-3.5973,3.5973]$ they $(2.55\times10^{-3},l_{\infty})$-coincide , while at interval $[-4.5973,4.5973]$ they $(5.97\times10^{-2},l_{\infty})$ - coincide.
The two relevant plots are given at Figure \ref{fig:01}.
\end{example}
\section{Practical computation of the effective radius of convergence}
The practical question is how can we find a radius such that the two graphs satisfactory coincide, despite the measure of their distance, at least for the beginning. In other words, if we just have a truncated series that we know has come from a smooth function, what is a first estimation for its useful range? We shall give a method based on the root test criterion for convergence.   
\subsection{Direct $n^{th}$-root method}
Since we know that for a convergent power series 
\begin{equation}
S(x)=\sum_{n=0}^{\infty}{a_n\,x^n} 
\end{equation}
the radius of convergence is given by
\begin{equation}
\label{eq:rdef}
R=\frac{1}{\limsup\limits_{n\rightarrow\infty}\sqrt[n]{\left|a_n\right|}}
\end{equation}
it is reasonable to argue that, approximately, when $n$ is not too small, the \ref{eq:rdef} has been reached if we take just the last term of the sequence
\begin{equation}
\label{eq:rndef}
\left\lbrace R_n=\frac{1}{\sqrt[n]{\left|a_n\right|}},\,\,n=1,2,\ldots,m\right\rbrace 
\end{equation}
But we can study further our series since for every function we can find an even and an odd function such that it can be written as the sum of those two parts.
\begin{equation}
\label{eq:feo}
\begin{matrix}
f_{ev}(x)&=&\frac{f(x)+f(-x)}{2} \\
f_{od}(x)&=&\frac{f(x)-f(-x)}{2} \\
f(x)&=&f_{ev}(x)+f_{od}(x)
\end{matrix}
\end{equation}
Now we can create two series, one with the even powers and another with odd powers. The same can be done for the truncated series, so for \ref{eq:pn} we have that
\begin{equation}
\begin{matrix}
\label{eq:seo}
P_{ev}(x)&=&a_0+a_2\,x^2+a_4\,x^4+\ldots \\
P_{od}(x)&=&a_1\,x+a_3\,x^3+a_5\,x^5+\ldots \\
P(x)&=&P_{ev}(x)+P_{od}(x)
\end{matrix}
\end{equation}
By this procedure we compute the relevant sequences for approximating the radius of convergence of both truncated series as
\begin{equation}
\label{eq:reo}
\begin{Bmatrix}
R_{2n}&=&\frac{1}{\sqrt[2n]{\left|a_{2n}\right|}} \\
&&\\
R_{2n+1}&=&\frac{1}{\sqrt[2n+1]{\left|a_{2n+1}\right|}} \\
\end{Bmatrix}   \,\,n=0,1,2,\ldots
\end{equation}
Finally we have two sets of values and if our approximation comes from a convergent series then both of them have to convergence to the same and true radius $R$. By the `convergence' here we actually mean the last term of every \ref{eq:reo} sequence, so our estimators are
\begin{equation}
\label{eq:reoest}
\begin{matrix}
\hat{R}_{ev}&=&R_{2k} \\
&&\\
\hat{R}_{od}&=&R_{2l+1}
\end{matrix}   
\end{equation}
where $k,l$ are defined from the last even and odd term respectively of the truncated series.

If the underlying function is even, then the sequence for the even truncated series has to be taken, otherwise if it is odd we examine the odd sequence. Another well known property of convergent power series is that their $n^{th}$ term has to vanish as $n\rightarrow\infty$, so we have a second check for the choice of the convergent sequence between the even and the odd one: we can disregard the series that comes from a divergent subset of coefficients.

\subsection{OLS regression $n^{th}$-root method}
If we observe that
$$
\frac{1}{R_n}=\sqrt[n]{\left|a_n\right|}={\left|a_n\right|}^{\frac{1}{n}}
$$
or
$$
\ln\left(\left|a_n\right|\right)=\left(-\ln\left(R_n\right)\right)\,n=\beta_{1}\,n
$$
we can obtain the $\hat{\beta}_{1}$ via OLS estimation for the data
\begin{equation}
\label{eq:datrr}
\left\{[n,\ln\left(\left|a_n\right|\right)] ,\,\,\forall a_n\neq{0}\right\}
\end{equation}
and then find an estimation for the radius of convergence as
\begin{equation}
\label{eq:rreg}
\hat{R}=e^{-\hat{\beta}_{1}}
\end{equation}
In this method sometimes it is better to first plot all data \ref{eq:datrr} in order to decide which subset of it has a better linear behaviour and then do OLS only for that.
\section{Numerical examples}
\begin{figure}
\begin{center}
\caption{Pdf of $N(0,1))$ and and its $10^{th}$ degree Taylor polynomial with $R_{ef}\left(0.1377,l_{\infty}\right)=2.3$}\label{fig:02}
\vspace{0.5in}  \includegraphics[width=6cm,height=6cm]{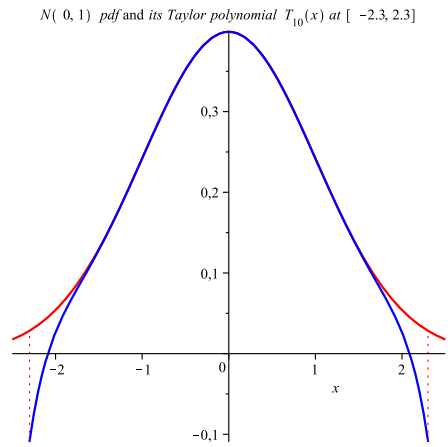}
\vspace{0.5in}
\end{center}
\end{figure}

\begin{figure}
\caption{Effective radius of convergence: $P_{10}(x)$ of $N(0,1)$ pdf}\label{fig:03} 
        \centering
        \begin{subfigure}[b]{0.4\textwidth}
                \includegraphics[width=\textwidth]{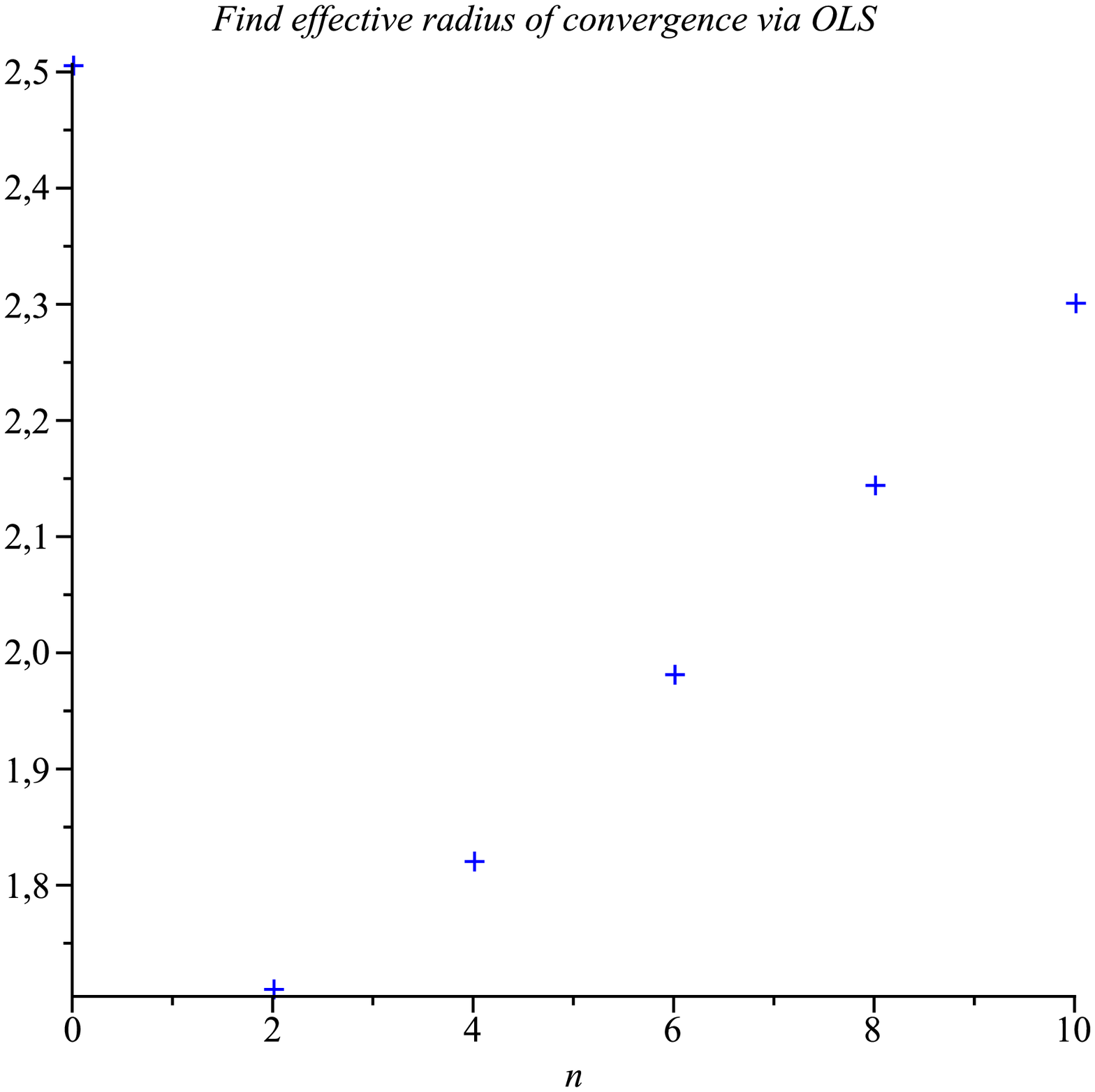}
                \caption{OLS data}
                \label{fig:03a}
        \end{subfigure}%
        \qquad
        \begin{subfigure}[b]{0.4\textwidth}
                \includegraphics[width=\textwidth]{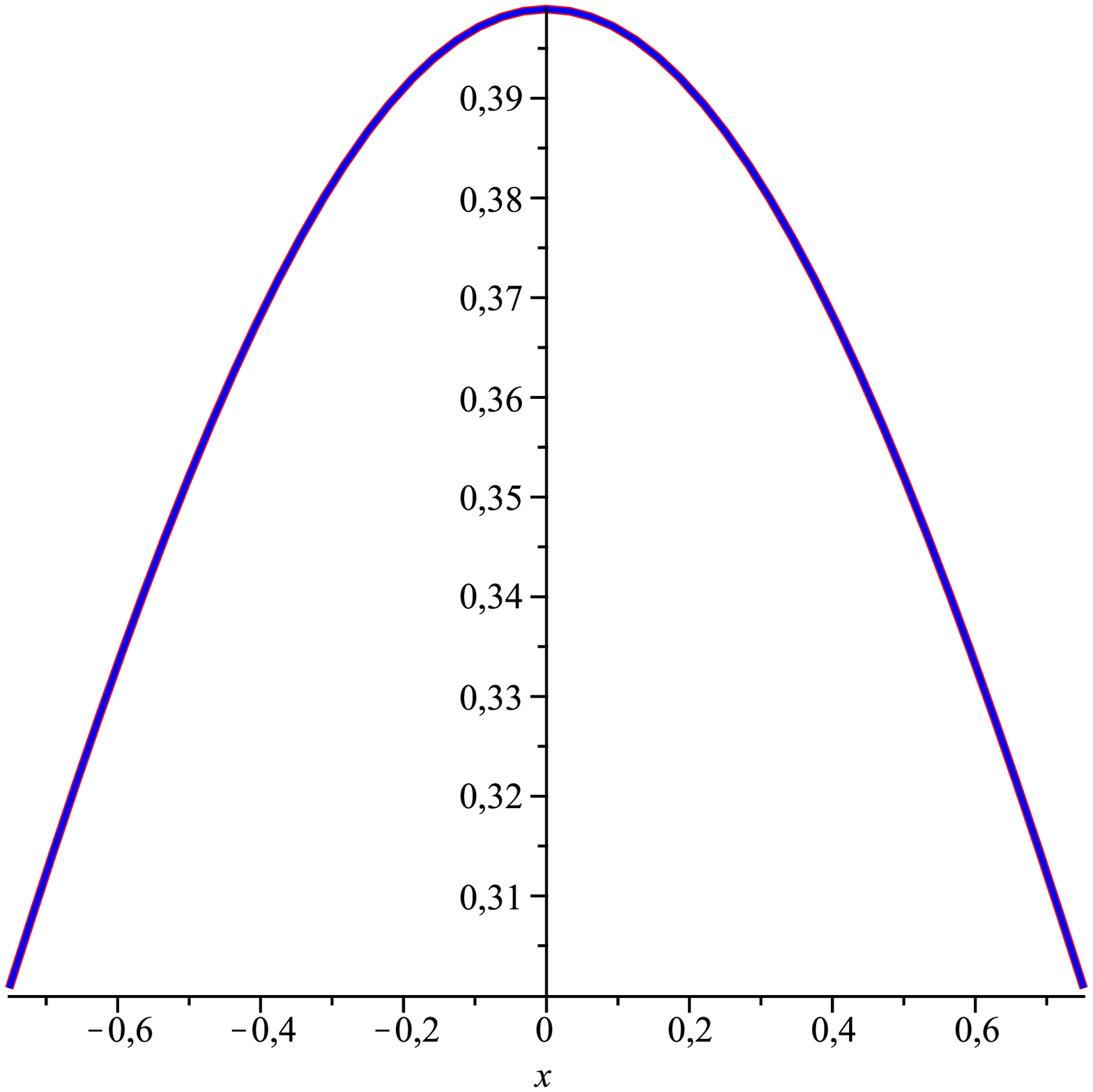}
                \caption{$R_{ef}(2.71\times10^{-7},l_{\infty})=0.752$}
                \label{fig:03b}
        \end{subfigure}
      
\end{figure}

\begin{figure}
\caption{Effective radius of convergence: $P_{11}(x)$ of $sin(x)$}\label{fig:04} 
        \centering
        \begin{subfigure}[b]{0.4\textwidth}
                \includegraphics[width=\textwidth]{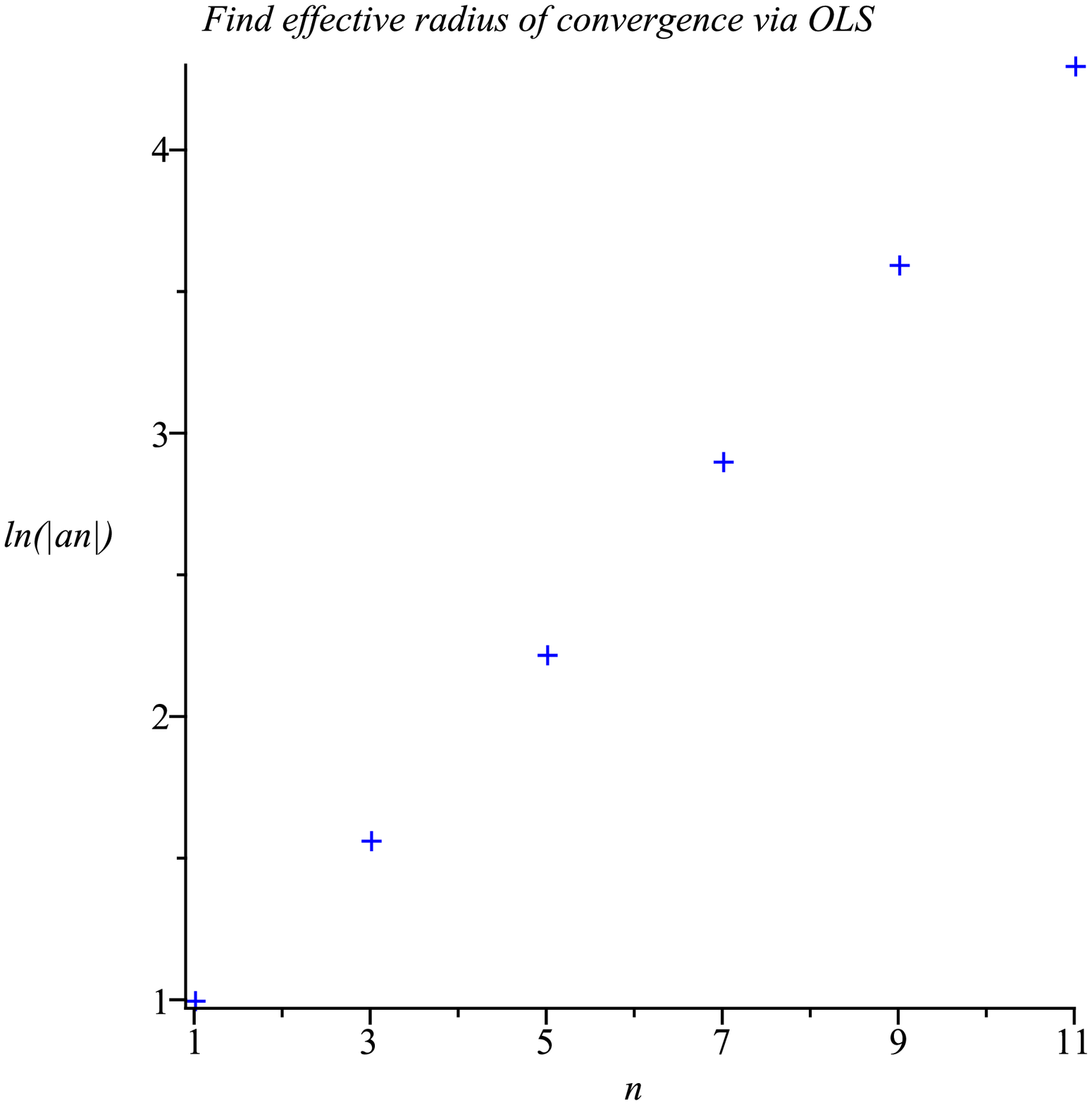}
                \caption{OLS data}
                \label{fig:04a}
        \end{subfigure}%
        \qquad
        \begin{subfigure}[b]{0.4\textwidth}
                \includegraphics[width=\textwidth]{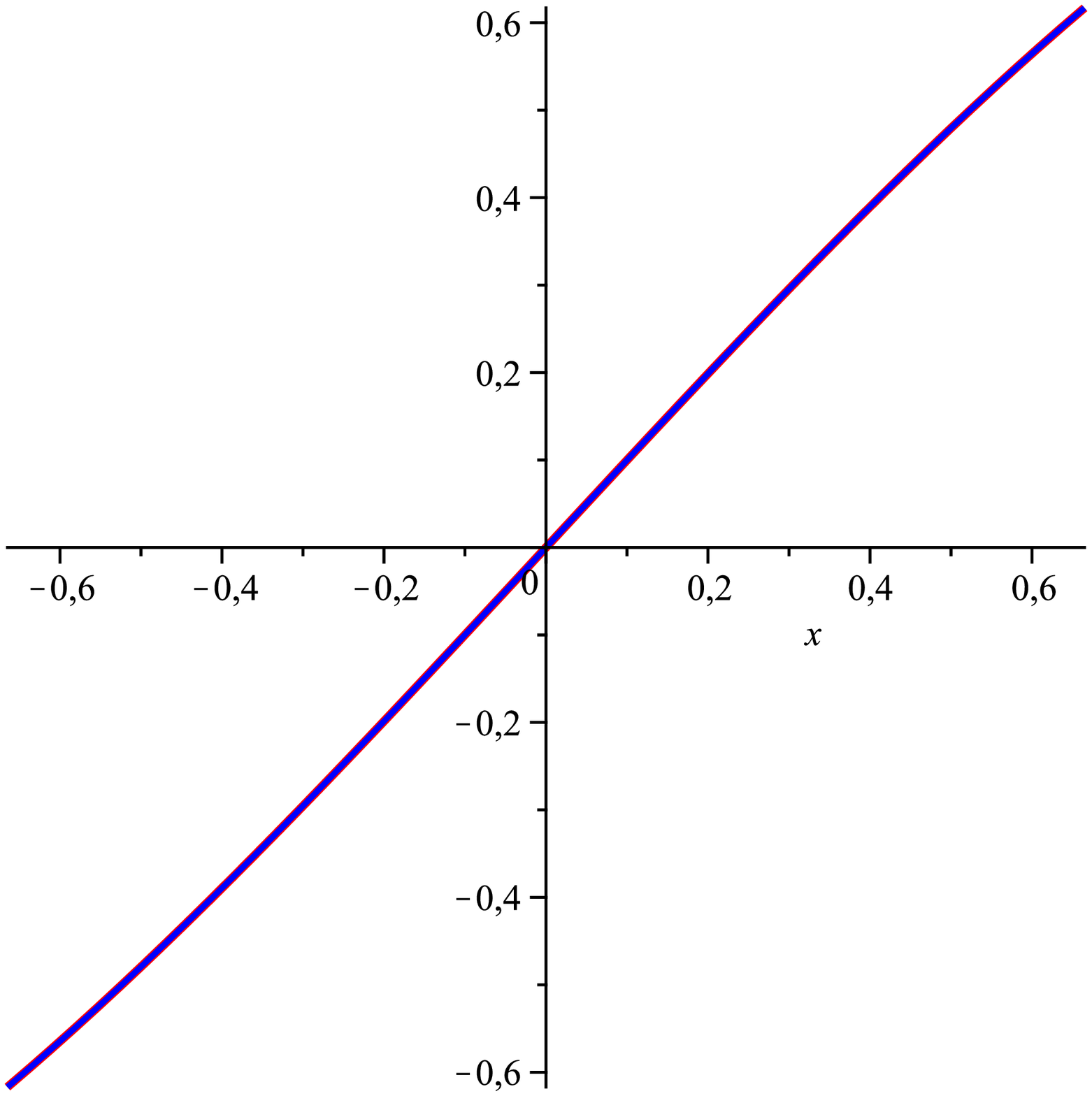}
                \caption{$R_{ef}(7.93\times10^{-13},l_{\infty})=0.665$}
                \label{fig:04b}
        \end{subfigure} 
\end{figure}

Let's illustrate our methods by studying a set of functions, even, odd or without parity.
\begin{example}
Let's examine the normal distribution pdf at interval $[-3,3]$ by both two methods. The function and its Taylor polynomial of degree 10 are
\begin{equation}
\label{eq:npdf}
\begin{array}{lll}
f \left( x \right)&=&\frac{1}{\sqrt{2\,\pi}}\,e^{-\frac{x^2}{2}}\\
&&\\
P_{12}(x)&=&
{\frac {1}{\sqrt {2\,\pi}}}-\frac{1}{2}\,{\frac {{x}^{2}}{\sqrt {
2\,\pi}}}+\frac{1}{8}\,{\frac {{x}^{4}}{\sqrt {2\,\pi}
}}-\frac{1}{48}\,{\frac {{x}^{6}}{\sqrt {2\,\pi}}}+{\frac {1}{384}
}\,{\frac {{x}^{8}}{\sqrt {2\,\pi}}}-{\frac {1}{3840}}\,{
\frac {{x}^{10}}{\sqrt {2\,\pi}}}
\end{array}
\end{equation}
\begin{itemize}
\item[1]
The sequence of even estimations is
$$
[ 2.507, 1.711, 1.822, 1.982, 2.145, 2.302]
$$
so we find $\hat{R}_{ev}=2.3=R_{ef}\left(0.1377,l_{\infty}\right)$ which is Figure \ref{fig:02}.
\item[2]
The OLS estimation is
$$
\hat{R}=0.7517
$$
The linear data used are shown at Figure \ref{fig:03a} while the plots for $x\in[-0.752,0.752]$ is Figure \ref{fig:03b}. \\
We can find that $\hat{R}_{ev}=0.752=R_{ef}\left(2.71\times 10^{-7},l_{\infty}\right)$, which gives a very small error, but our useful range has decreased too much.
\end{itemize}
\end{example}
%
%
\begin{example}
Let's examine the sinus function with both methods.
\begin{itemize}
\item[1]
Our odd sinus function \ref{eq:sin} and its truncated series \ref{eq:sin11}  gives as the sequence of odd estimations
$$
[ 1.0, 1.565, 2.221, 2.903, 3.597, 4.300]
$$
so we directly find $\hat{R}_{od}=4.3=R_{ef}\left(2.53\times 10^{-2},l_{\infty}\right)$ which is just what we have already observed at Figure \ref{fig:00}.
\item[2]
The OLS estimation is
$$
\hat{R}=0.665
$$
which came from an almost perfect straight line shown at Figure \ref{fig:04a} while the plots for the interval $[-0.665,0.665]$ is Figure \ref{fig:04b}. For this approximation we have that $\hat{R}_{od}=0.665=R_{ef}\left(7.93\times 10^{-13},l_{\infty}\right)$ which again gives a very small error as compensating us for the smaller useful range.
\end{itemize}
\end{example}

\begin{figure}
\begin{center}
\caption{$f \left( x \right) =\sin \left( 3\,x \right) \cos \left( 5\,x \right)
{e^{-x}}+3\,\sin \left( \pi \,x \right) {e^{\frac{x}{2}}}$  and its $30^{th}$ degree Taylor polynomial with $R_{ef}\left(0.1377,l_{\infty}\right)=1.54$}\label{fig:05}
\vspace{0.5in}  \includegraphics[width=6cm,height=6cm]{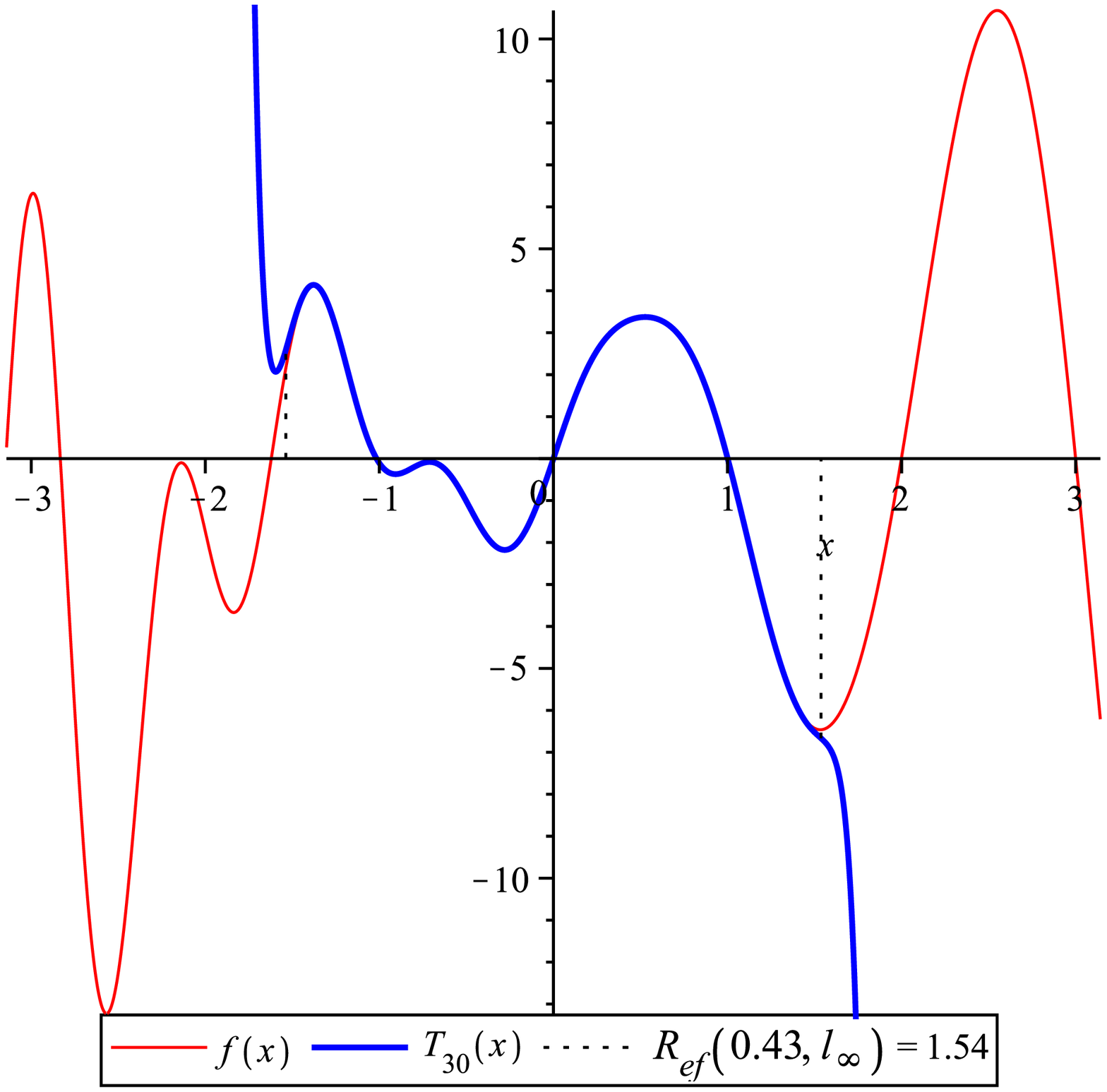}
\vspace{0.5in}
\end{center}
\end{figure}

\begin{figure}
\caption{$f \left( x \right) ={\frac {\frac{1}{8}\,x+\frac{1}{2}\,{x}^{2}}{1+\frac{1}{8}\,x+\frac{1}{2}\,{x}^{2}}}$ and its $30^{th}$ degree Taylor polynomial with effective radius of convergence $R_{ef}\left(0.217,l_{\infty}\right)=1.402$}\label{fig:06} 
        \centering
        \begin{subfigure}[b]{0.4\textwidth}
                \includegraphics[width=\textwidth]{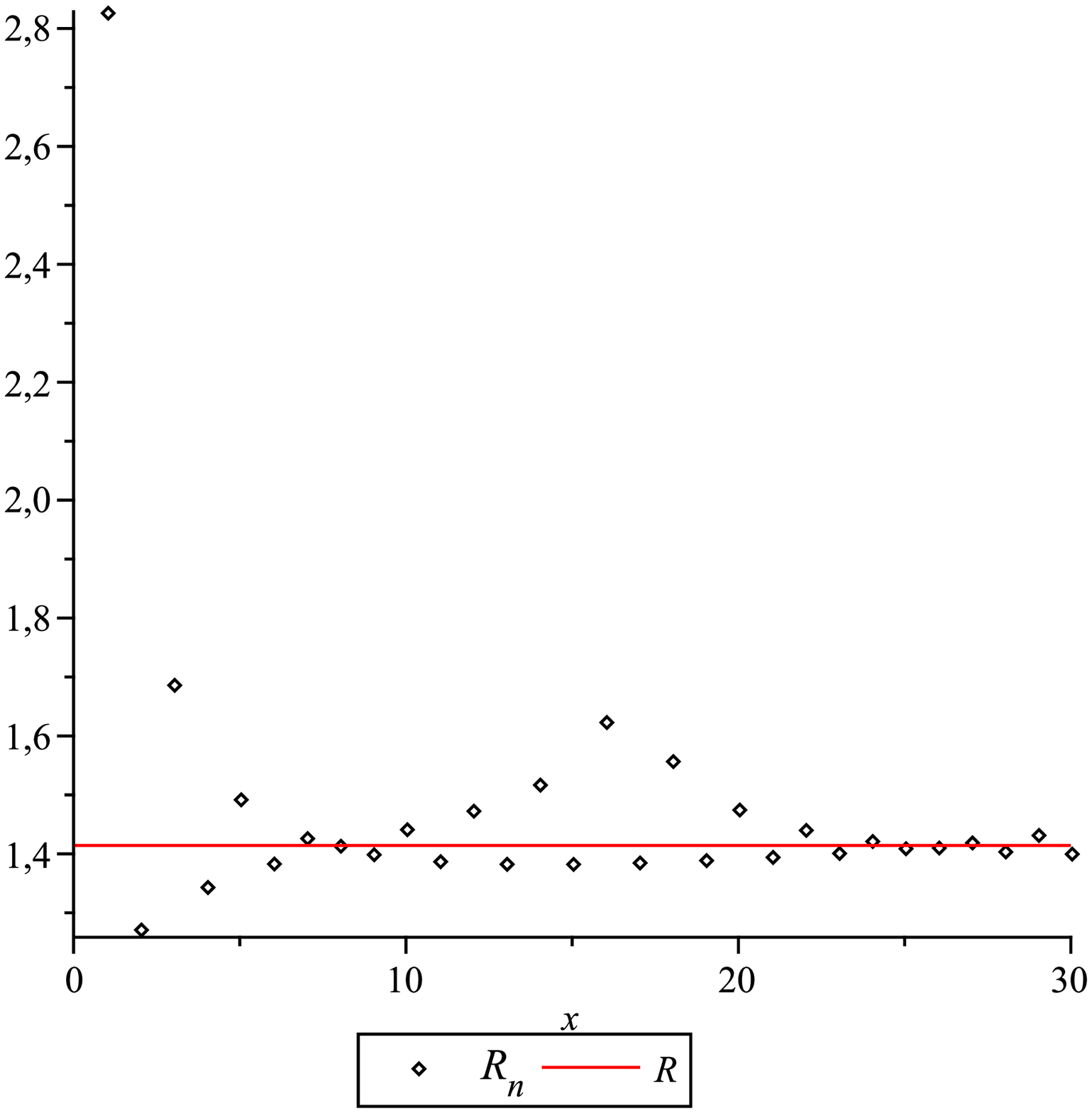}
                \caption{$R_n\rightarrow{R}$}
                \label{fig:06a}
        \end{subfigure}%
        \qquad
        \begin{subfigure}[b]{0.4\textwidth}
                \includegraphics[width=\textwidth]{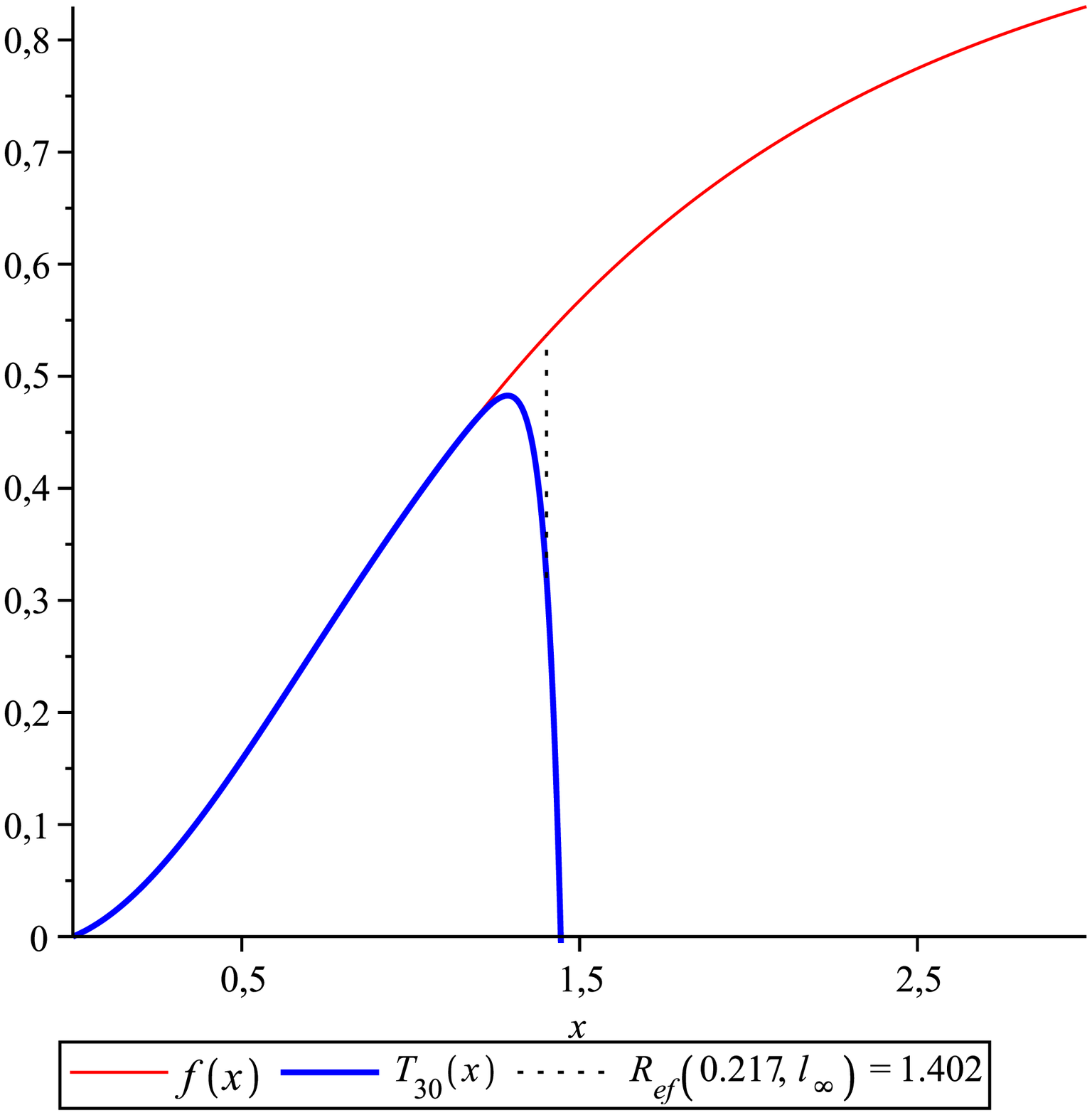}
                \caption{$R_{ef}(0.217,l_{\infty})=1.402$}
                \label{fig:06b}
        \end{subfigure}
        \qquad
                \begin{subfigure}[b]{0.4\textwidth}
                        \includegraphics[width=\textwidth]{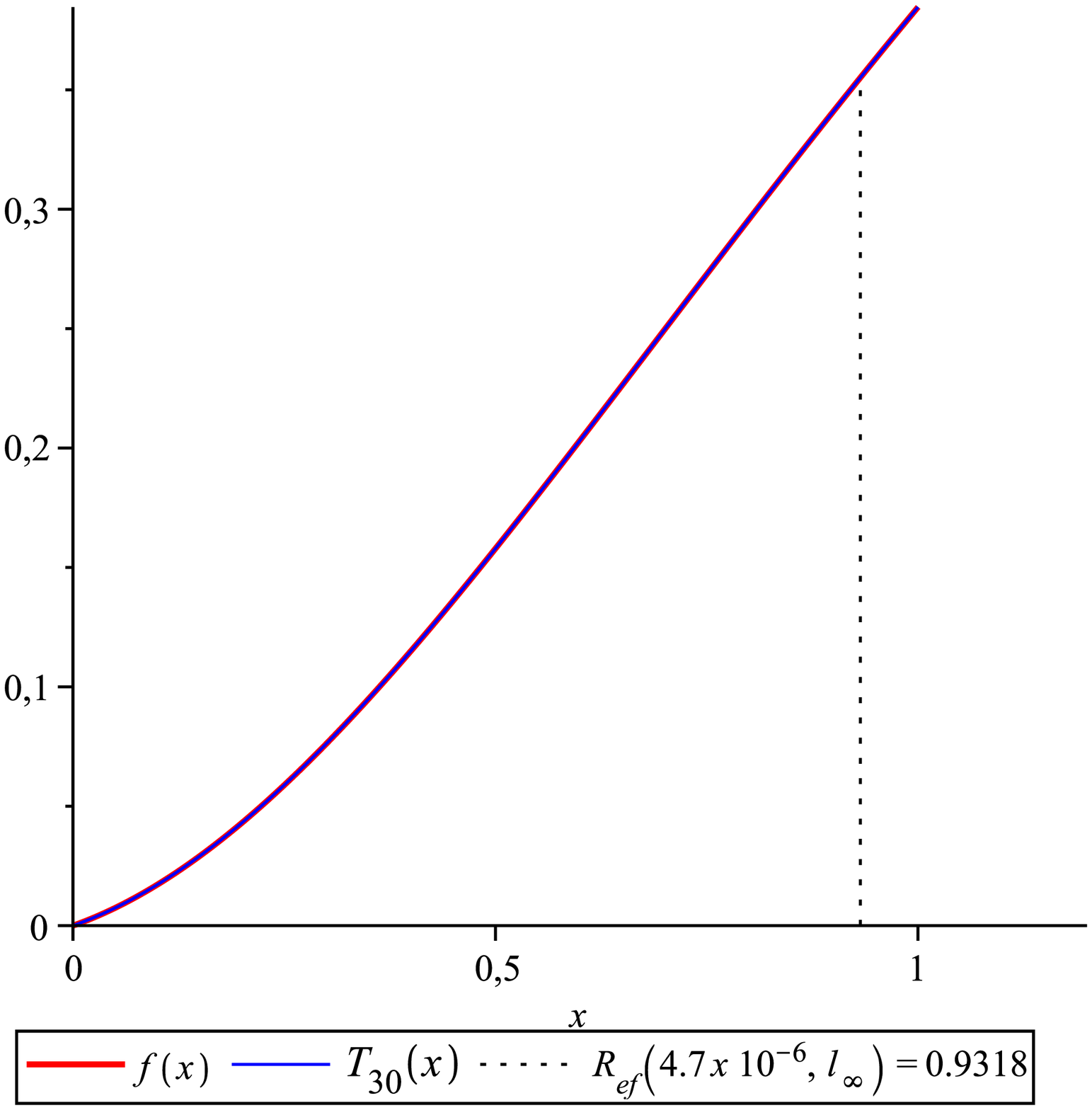}
                        \caption{$R_{ef}(4.74\times10^{-7},l_{\infty})=0.9318$}
                        \label{fig:06c}
                \end{subfigure}
      
\end{figure}

\begin{example}
Let's study now a more complicated function which is not even or odd, see \cite{dch-13} for more details. We find that the effective radius of convergence is $R_{ef}\left(0.1377,l_{\infty}\right)=1.54$ as clearly is shown at Figure \ref{fig:05}.\\ We can also compute the relevant OLS estimation and find that is  $\hat{R}=0.946=R_{ef}\left(1.24 \times 10^{-7},l_{\infty}\right)$, which is less useful since it is only approximately the $\frac{2}{3}$ of the root test based one. 
\end{example}
An interesting note has to be done here about the relationship between the real radius of convergence and the effective radius that we compute via the sequences \ref{eq:rndef} and \ref{eq:reo}. 
\begin{example}
As an example we shall try to find the $R_{ef}$ of a $2:2$-function according to  \cite{bar-75} notation:
\begin{equation}
\label{eq:22}
f \left( x \right) ={\frac {\frac{1}{8}\,x+\frac{1}{2}\,{x}^{2}}{1+\frac{1}{8}\,x+\frac{1}{2}\,{x}^{2}}}
\end{equation}
If we take its $30^{th}$ degree Taylor polynomial and compute our sequence \ref{eq:rndef} we see that it convergences to the true radius of convergence
 $$
 R=\left| -\frac{1}{8}\pm i\,\frac{1}{8}\,\sqrt {127}\right| =1.41421356
 $$
as is presented at Figure \ref{fig:06a} while $R_{ef}\left(0.217,l_{\infty}\right)=1.40198948$ is shown at Figure \ref{fig:06b}. The OLS estimation is $\hat{R}=0.9318=R_{ef}\left(4.74\times 10^{-7},l_{\infty}\right)$ and for the interval $[0,0.9318]$ the function and its Taylor approximation are indistinguishable, see Figure \ref{fig:06c}.
\end{example}
\section{Discussion}
We saw that for every truncated Taylor power series we can find two radius, one based directly on a root test sequence and another using a proper OLS linear regression on those values. The first one always gives the greatest range $[x_{0}-R_{ef},  x_{0}+R_{ef}]$ where we can plot our function $f(x)$ and its Taylor polynomial $T_m(x)$ with almost identical curves. The second gives a more restrictive range, but with very small error measured by $l_{\infty}$ norm.\\
We have also found that our effective radius of convergence, if we take the degree of the Taylor polynomial sufficiently large, is a rough estimation of the true non-infinite radius of convergence.

\newpage

\begin{thebibliography}{2}

\bibitem{dch-13}{\sc D.T. Christopoulos},
{Polynomial regression using trapezoidal rule for computing Legendre coefficients, {http://arxiv.org/abs/1311.7525v1}, 2013 }

\bibitem{bar-75}
{\sc W. Bardsley \& R. Childs}, {Sigmoid Curves, Nonlinear Double-Reciprocal Plots and Allosterism, \emph{Biochemical Journal}, \textbf{149}, 313-328, 1975}

\end{thebibliography}
\end{document}